\documentclass[12pt]{article}
\usepackage{mathrsfs}
\usepackage{amsfonts}
\usepackage{amssymb,amsmath}
\usepackage{cite}
\usepackage{cases}
\def\cl{\centerline}
\def\la{\lambda}
\def\al{\alpha}
\def\vs{\vspace*}

\def\Z{\mathbb{Z}}
\def\C{\mathbb{C}}
\def\QED{\hfill$\Box$}
\def\ni{\noindent}

\numberwithin{equation}{section}
\newtheorem{theo}{Theorem}[section]

\newtheorem{lemm}[theo]{Lemma}

\newtheorem{case}{Case}
\newtheorem{step}{Step}
\newtheorem{subcase}{Subcase}

\begin{document}
\cl{{\bf\large Derivations and automorphisms of twisted}}
\vs{6pt}
\cl{{\bf\large  deformative
Schr\"{o}dinger-Virasoro Lie algebras}\footnote
{\!\!Supported by NSF grants 10825101, 10861004 of China, China Postdoctoral Science Foundation grant 201003326 and the Natural Science Research Projects of Jiangsu Education Committee grant 09KJB110001\\[2pt]
\indent\ \,Electronic mail: wwll@mail.ustc.edu.cn, sd\_junbo@163.com, xying@mail.ustc.edu.cn}}

\vs{12pt}

\cl{Wei Wang$^{1,2)}$,  Junbo Li$^{2,3)}$, Ying Xu$^{2)}$}

\cl{\footnotesize $^{1)}$School of Mathematics and Computer Science,
Ningxia University, Yinchuan 750021, China}

\cl{\footnotesize$^{2)}$Wu Wen-Tsun Key Laboratory of Mathematics,
USTC, Hefei 230026, China}

\cl{\footnotesize $^{3)}$School of Mathematics and Statistics,
Changshu Institute of Technology, Changshu 215500, China}

\vs{8pt}

{\small
\parskip .005 truein
\baselineskip 3pt \lineskip 3pt

\noindent{{\bf Abstract:}  In this paper the derivation algebra and
automorphism group of the twisted deformative
Schr\"{o}dinger-Virasoro Lie algebras are determined.
\vs{5pt}

\ni{\bf Key words:} Schr\"{o}dinger-Virasoro Lie algebras, derivations, automorphisms.}

\ni{\it Mathematics Subject Classification (2010):} 17B05, 17B40, 17B65, 17B68.}
\parskip .001 truein\baselineskip 6pt \lineskip 6pt
\section{Introduction}
\setcounter{section}{1}\setcounter{equation}{0}

\indent\ \ \ \ \ It is well known that the Virasoro algebra  plays an important role in many areas of theoretical physics and mathematics, which occures in the investigation of conformal field theory and has a $\C$-basis $\{L_n,c\,|\,n\in\Z\}$ with the nontrivial relations $[L_n,L_m]=(m-n)L_{n+m}$. It can be regarded as the complexification of the Lie algebra of polynomial vector fields on a circle, and also as the Lie algebra of derivations of the ring $\C[z,z^{-1}]$. The centerless Virasoro algebra admits many kinds of extensions, one of these is the Schr\"{o}dinger-Virasoro type Lie algebras (see \cite{MJ1,MJ2,RU}), firstly introduced in \cite{H} in the context of non-equilibrium
statistical physics during the process of investigating the free Schr\"{o}dinger equations and closely related to the Schr\"{o}dinger algebra and the Virasoro algebra. Recently the vertex algebra representations, Lie bialgebra structures, irreducible weight modules with finite-dimensional weight spaces and also Wittaker modules of the Schr\"{o}dinger-Witt Lie algebras were extensively investigated in \cite{U,HLS,LS1,ZT1}. The generalization of the Schr\"{o}dinger-Virasoro Lie algebras was introduced in \cite{TZ1}, whose automorphism groups and Verma modules were described therein.

For any $\la$, $\mu\in\C$, \cite{RU} introduced a family of  infinite-dimensional Lie algebras called {\it  twisted deformative Schr\"{o}dinger-Virasoro Lie algebras}, admitting $\C$-basis
$\{L_n,Y_n,M_n\,|\,n\in\Z\}$ and the following Lie brackets
\begin{eqnarray*}
&&[L_n,L_m]=(m-n)L_{n+m},\\
&&[L_n,Y_m]=(m-\frac{\la+1}{2}n+\mu)Y_{n+m},\ \ \ [Y_n,Y_m]=(m-n)M_{n+m},\\
&&[L_n,M_m]=(m-\la n+2\mu)M_{n+m},\ \ \ \ \ [Y_n,M_m]=[M_n,M_m]=0.
\end{eqnarray*}
We denote this Lie algebra by $\mathscr{L}$, which is $\Z$-graded with
\begin{equation*}
\mathscr{L}=\bigoplus_{n\in\Z}\C\mathscr{L}_n,\ \ \mathscr{L}_n=\C L_n\oplus\C Y_n\oplus \C M_n,\ \ [\mathscr{L}_n,\mathscr{L}_m]\subseteq\mathscr{L}_{n+m}.
\end{equation*}
For convenience we introduce the following notations
\begin{eqnarray*}
\mathcal{L}=\mbox{$\sum\limits_{n\in\Z}$}\C L_n,\ \ \ \mathcal{Y}=\mbox{$\sum\limits_{n\in\Z}$}\C Y_n,\ \ \mathcal{M}=\mbox{$\sum\limits_{n\in\Z}$}\C M_n.
\end{eqnarray*}
Then $\mathcal{L}$  is the centerless Virasoro algebra, $\mathcal{I}=\mathcal{Y}\oplus \mathcal{M}$ is the unique maximal ideal of $\mathscr{L}$ and $\mathcal{M}$ is the center of $\mathcal{I}$.

It is well known that the determination of derivations and
automorphisms plays an important part in the investigation of the structure and representation of the relevant Lie algebras. Many references (see \cite{F1,DZ,RU,GJP1,S1,S2,S3,SJ}) have focused on derivations and automorphisms of different Lie algebra backgrounds.

Note that if $\mu\in\Z$, then $\{L_n,Y_{n-\mu},M_{n-2\mu}\}$ is a basis of $\mathscr{L}$. Hence one can assume $\mu=0$ if $\mu\in\Z$. Especially, for the case $\la=\mu=0$, $\mathscr{L}$ is nothing but the twisted Schr\"{o}dinger-Virasoro algebra, whose  derivations and automorphisms were determined in \cite{LS2}. If $\mu\in \frac{1}{2}+\Z$, one can assume $\mu=\frac{1}{2}$ by shifting basis, in which case the derivations were investigated in \cite{ZH} and automorphisms for the special case $\la=0$ were determined in \cite{TZ1}. Furthermore, for the case $\mu=\frac{1}{2}$ and $\la\neq0$, one can obtain the corresponding results on automorphisms by following the proof of Theorem 3.2 in \cite{TZ1}. Thus, in this paper we always make the following assumptions on $\mu$ and $\la$
\begin{equation}\label{mu}
\mu\notin\frac{1}{2}+\Z,\ \ \mu=0\ \ {\rm and}\ \ \la\neq0\ \ {\rm if}\ \ \mu\in\Z.
\end{equation}

In the following two sections, we shall determine the derivation algebra (see Theorem \ref{theorem1}) and the automorphism group (see Theorem \ref{theo}) of $\mathscr{L}$ under the assumptions made in (\ref{mu}).

\section{ Derivation algebra of $\mathscr{L}$}
\setcounter{section}{2}\setcounter{theo}{0}

\indent\ \ \ \ \ A linear map $d:\mathscr{L}\rightarrow \mathscr{L}$ is called a {\it derivation} of $\mathscr{L}$, if $d\big([x,y]\big)=[d(x),y]+[x, d(y)]$ holds for any $x$, $y\in\mathscr{L}$. For any fixed $z\in\mathscr{L}$, the linear map $ad_{z}\!:\mathscr{L}\rightarrow\mathscr{L}$ is called an {\it inner} derivation if $ad_{z}(x)=[z,x]$ for any $x\in\mathscr{L}$. Denote by \emph{Der\,$\mathscr{L}$} and \emph{ad$\,\mathscr{L}$} respectively the vector spaces of all derivations and inner derivations. Then the first homological group
$\mathcal{H}^{1}(\mathscr{L},\mathscr{L})\cong\emph{Der}\mathscr{L}/{ad}\,\mathscr{L}$.

Firstly, we give the following description of \emph{Der\,$\mathscr{L}$}.
\begin{lemm}\ \ \label{lemma1}
Der$\hspace{0.08cm}\mathscr{L}=(Der\mathscr{L})_0+ad\,\mathscr{L}$,
where
\begin{eqnarray*}
&&(Der\mathscr{L})_0=\{d\in\mbox{Der}\hspace{0.08cm}\mathscr{L}\
|\ d(\mathscr{L}_n)\subseteq\mathscr{L}_n,\,\forall\,n\in\Z\}.
\end{eqnarray*}
\end{lemm}
{\bf Proof}\ \ \ By Proposition $1.2$ of \cite{F1},  we need two steps to complete the proof of this lemma.
\begin{step}
For any $n\neq0$, if $d\in\mathcal{H}^{1}(\mathscr{L}_0,\mathscr{L}_n)$, then $d\in{ad}{\,\mathscr{L}}$.
\end{step}
Since $\mathscr{L}_n=\C L_n\oplus\C Y_n\oplus\C M_n$ for all $n\in\Z$, one can assume $d(X_0)=e^X_1L_n+e^X_2Y_n+e^X_3M_n$
for some $e^X_i\in\mathbb{C}$, $i=1,2,3$, $X\in\{L,Y,M\}$. Applying $d$ to $[L_0,Y_0]=\mu Y_0$, $[L_0,M_0]=2\mu M_0$ and $[Y_0,M_0]=0$ respectively,  we obtain
\begin{eqnarray*}
\left\{
\begin{array}{llll}
(n-\mu) e_1^Y=(\mu-\frac{1+\la}{2}n)e_1^L+ne_2^Y=ne_2^L -(n+\mu)e_3^Y=0,\vs{6pt}\\
(n-2\mu)e_1^M=(n-\mu)e_2^M=(2\mu-\la n)e_1^L- ne_3^M=0,\vs{6pt}\\
(2\mu-\la n)e_1^Y+ne_2^M=(\mu-\frac{1+\la}{2}n)e_1^M=0,
\end{array}\right.
\end{eqnarray*}
which combined with (\ref{mu}) give
\begin{eqnarray*}
&& e_1^Y=e_1^M=e_2^M=0,\ \ e_2^Y=\displaystyle\frac{\la-2\mu+1}{2n}e_1^L,\ \ e_3^Y=\displaystyle\frac{n}{n+\mu}e_2^L,\ \ e_3^M=\displaystyle\frac{2\mu-\la n}{n}e_1^L.
\end{eqnarray*}
Hence, denoting $\al=\frac{1}{n}e_1^LL_n+\frac{1}{n+\mu}e_2^LY_n+\frac{1}{n+2\mu}e_3^LM_n$, we obtain $d(X_0)=[X_0,\al]$ for any $X\in\{L,Y,M\}$.
\begin{step}
For any $n\neq m$, if $f\in{\rm Hom}_{\mathscr{L}_{0}}(\mathscr{L}_n,\mathscr{L}_m)$, then $f=0$.
\end{step}
Assume that $f(X_n)=c^X_1L_m+c^X_2Y_m+c^X_3M_m$ for some $c^X_i\in\C$ with $X\in\{L,Y,M\}$, $i=1,2,3$.
Applying $f$ to both sides of the following three identities
\begin{eqnarray*}
[L_0,L_n]=n L_n,\ [L_0,Y_n]=(n+\mu)Y_n,\  [L_0,M_n]=(n+2\mu)M_n,
\end{eqnarray*}
and comparing the coefficients of $L_m$, $Y_m$ and $M_m$ respectively, we have
\begin{eqnarray*}
\left\{
\begin{array}{llll}
(m-n)c_1^L=(m-n+\mu)c_2^L=(m-n+2\mu)c_3^L=0,\vs{4pt}\\
(m-n-\mu)c_1^Y=(m-n)c_2^Y=(m-n+\mu)c_3^Y=0,\vs{4pt}\\
(m-n-2\mu)c_1^M=(m-n-\mu)c_2^M=(n+2\mu)c_3^M =0.
\end{array}\right.
\end{eqnarray*}
Note that $\mu$ satisfies (\ref{mu}), one can deduce $c^X_i=0$, $X\in\{L,Y,M\}$, $i=1,2,3$. Thus $f=0$.\QED

Let $d\in (Der\mathscr{L})_0$. For $n\in\Z$ and $X\in\{L,Y,M\}$,  one can assume
\begin{equation*}
d(X_n)=f^X_1(n)L_n+f^X_2(n)Y_n+f^X_3(n)M_n\ \ \ \ {\rm for\ \,some}\ f^X_i(n)\in\C,\ i=1,2,3.
\end{equation*}
\begin{lemm}\label{lemma2}   \ \ For any $n\in\Z$ and some $a$, $b$, $c$, $\bar{c}$, $e$,  $\bar{e}\in\C$, one can assume

(1)\ \ $f_1^L(n)=an$.

(2)\ \ $f_2^L(n)=\left\{\begin{array}{llllllll}
\frac{e}{2\mu}\big(2\mu-(\la+1)n\big)&{\rm if}\ \ \mu\notin\Z,\vs{6pt}\\
bn(1-\delta_{\la,-1})&{\rm if}\ \ \mu=0.
\end{array}\right.$

(3)\ \ If $\mu\notin\Z$, then $f_3^L(n)=\frac{\bar{e}}{2\mu}(2\mu-\la n)$.

(4)\ \ If $\mu=0$, then

\ \ \ \ \ \ \ \ \ \ \ \ \ $f_3^L(n)=\left\{\begin{array}{llllllll}
\frac{c}{6}(n^3-n)-\frac{\bar{c}}{3} (n^3-4n)&{\rm if}\ \ \la=-2,\vs{6pt}\\
\frac{c}{2}(n^2-n)-\bar{c}(n^2-2n)&{\rm if}\ \ \la=-1,\vs{6pt}\\
\bar{c}n&{\rm if}\ \ \la\notin\{-2,0,-1\}.
\end{array}\right.$
\end{lemm}
\ni{\bf Proof}\ \ (1)\ \ Applying $d$ to $[L_n,L_m]=(m-n)L_{m+n}$ and comparing the coefficients of $L_{m+n}$, we have $(m-n)\big(f_1^L(m+n)-f_1^L(m)-f_1^L(n)\big)=0$, which gives \begin{equation}\label{44}
f_1^L(m+n)=f_1^L(m)+f_1^L(n)\ \ {\rm if}\ \ m\neq n.
\end{equation}
Taking $n=0$ in (\ref{44}), we have $f_1^L(0)=0$. Applying (\ref{44}), one can deduce
\begin{equation*}\label{43}
f_1^L(2n)=f_1^L(3n)+f_1^L(-n)=f_1^L(2n)+f_1^L(n)+f_1^L(-2n)+f_1^L(n)
=2f_1^L(n),
\end{equation*}
which together with (\ref{44}) gives $f_1^L(m+n)=f_1^L(m)+f_1^L(n)$ for all $m$, $n\in\Z$. Thus $f_1^L(n)=f_1^L(n-1)+f_1^L(1)$. Hence, by induction on $n$, one can deduce $f_1^L(n)=nf_1^L(1)$ for all $n\in\Z$.

(2)\ \ Applying $d$ to $[L_n,L_m]=(m-n)L_{m+n}$ and comparing the coefficients of $Y_{m+n}$, one has
\begin{eqnarray}
\label{2}\big(2m-(\la+1)n+2\mu)f^L_2(m)-\big(2n-(\la+1)m+2\mu\big)f^L_2(n)=2(m-n)f^L_2(m+n).
\end{eqnarray}
\begin{case}
$\mu\notin\Z$.
\end{case}
Setting $n=0$ in (\ref{2}), we obtain
$2\mu f^L_2(m)=\big(2\mu-(\la+1)m\big)f^L_2(0)$, which gives $f^L_2(m)=\frac{1}{2\mu}\big(2\mu-(\la+1)m\big)f^L_2(0)$ for all $m\in\Z$.

\begin{case}
$\mu=0$.
\end{case}
Taking $\mu=0$ in (\ref{2}), we have
\begin{eqnarray}
\label{8}\big(2m-(\la+1)n\big)f^L_2(m)-\big(2n-(\la+1)m\big)f^L_2(n)=2(m-n)f^L_2(m+n).
\end{eqnarray}
Taking $n=1,2$ in (\ref{8}), we obtain
\begin{eqnarray}
&&\label{12}2(m-1)f^L_2(m+1)=(2m-\la-1)f^L_2(m)-\big(2-(\la+1)m\big)f^L_2(1),\\
&&\label{40}2(m-2)f^L_2(m+2)=2(m-\la-1)f^L_2(m)-\big(4-(\la+1)m\big)f^L_2(2).
\end{eqnarray}
Taking $n=1$ and replacing $m$ by $m+1$ in (\ref{8}), one can deduce
\begin{eqnarray}\label{41}
(2m-\la+1)f^L_2(m+1)-\big(2-(\la+1)(m+1)\big)f^L_2(1)=2mf^L_2(m+2).
\end{eqnarray}
Multiplying (\ref{41}) by $2(m-1)(m-2)$, then using (\ref{12}) and (\ref{40}), we obtain
\begin{eqnarray}\label{42}
&&\ \ \ \nonumber(\la-1)\big(2\la - (5+ \la)m+2\big)f^L_2(m)\\
&&=m(m-2 ) \big(4(1 + \la) m- \la^2-7\big)f^L_2(1)+2m(m-1)\big(4-(\la+1)m\big)f^L_2(2).
\end{eqnarray}

Applying $d$ to $[L_n,Y_m]=(m-\frac{\la+1}{2}n)Y_{m+n}$ and comparing the coefficients of $M_{m+n}$, we have
\begin{eqnarray}
\label{Y3}2(m-n)f^L_2(n)+2(m-\la n)f^Y_3(m)-(2m-(1+\la)n)f^Y_3(n+m)=0.
\end{eqnarray}

Setting $m=0$ in (\ref{Y3}), we obtain
\begin{equation}\label{05}
(1+\la)nf^Y_3(n)=2nf^L_2(n)+2\la nf^Y_3(0).
\end{equation}
Replacing $n$ by $-n$, $m$ by $n$ in (\ref{Y3}), one has
\begin{eqnarray}\label{06}
(\la+3)nf^Y_3(0)=2(1+\la)nf^Y_3(n)+4nf^L_2(-n).
\end{eqnarray}
Multiplying (\ref{05}) by $(\la+3)$ and using (\ref{06}), we obtain
\begin{equation}\label{07}
3(1-\la^2)nf^Y_3(n)=2(\la+3)nf^L_2(n)+8\la nf^L_2(-n).
\end{equation}
Taking $n=1,2$ in (\ref{07}) respectively, one has
\begin{eqnarray}
&&\label{08}3(1-\la^2)f^Y_3(1)=2(\la+3)f^L_2(1)+8\la f^L_2(-1),\\
&&\label{09}3(1-\la^2)f^Y_3(2)=2(\la+3)f^L_2(2)+8\la f^L_2(-2).
\end{eqnarray}
Furthermore, setting $n=m=1$ in (\ref{Y3}), we obtain
\begin{equation}\label{L301}
(1-\la)\big(f^Y_3(2)-2f^Y_3(1)\big)=0.
\end{equation}
Multiplying (\ref{L301}) by $3(1+\la)$, using (\ref{08}) and (\ref{09}), one can deduce
\begin{equation}\label{main}
(\la+3)\big(f^L_2(2)-2f^L_2(1)\big)+4\la\big(f^L_2(-2)-2f^L_2(-1)\big)=0.
\end{equation}

\setcounter{subcase}{0}
\begin{subcase}
$\la=-5$.
\end{subcase}
Taking $\la=-5$ in (\ref{42}), one has
\begin{eqnarray}\label{010}
f^L_2(m)=\frac{1}{6}m(m-1)(m+1)f^L_2(2)-\frac{1}{3}m(m-2 )(m+2)f^L_2(1).
\end{eqnarray}
Taking $m=-1,-2$ in (\ref{010}), we obtain $f^L_2(-1)=-f^L_2(1)$ and $f^L_2(-2)=-f^L_2(2)$. Then (\ref{main}) gives $f^L_2(2)=2f^L_2(1)$, which together with (\ref{010}), admits $f^L_2(m)=mf^L_2(1)$ for all $m\in\Z$.
\begin{subcase}
$\la=-3$.
\end{subcase}
Taking $\la=-3$ in (\ref{42}), one has
\begin{equation}\label{21}
2(m+2)f^L_2(m)=m(m+2)(m-1)f^L_2(2)-2 m (m+2)(m-2)f^L_2(1),
\end{equation}
which gives
\begin{equation}\label{27}
f^L_2(m)=\frac{1}{2}m(m-1)f^L_2(2)-m(m-2)f^L_2(1)\ \ {\rm if}\ \ m\neq-2.
\end{equation}
Setting $\la=-3$, $m=1$ and $n=-2$ in (\ref{8}), we have
\begin{equation}\label{54}
f^L_2(-2)=f^L_2(1)+3f^L_2(-1).
\end{equation}
Furthermore, taking $m=-1$ in (\ref{27}), one has
\begin{equation}\label{003}
f^L_2(-1)=f^L_2(2)-3f^L_2(1).
\end{equation}
Applying this by substituting for the second term of the right-hand side in (\ref{54}), we have
\begin{eqnarray}\label{002}
f^L_2(-2)=3f^L_2(2)-8f^L_2(1).
\end{eqnarray}
Thus (\ref{27}) holds for all $m\in\Z$. Taking $\la=-3$ in (\ref{main}), then using (\ref{003}) and (\ref{002}), one can deduce $f^L_2(2)=2f^L_2(1)$. Thus (\ref{27}) gives  $f^L_2(m)=mf^L_2(1)$ for all $m\in\Z$.
\begin{subcase}
$\la=-1$.
\end{subcase}
Letting $\la=-1$ in (\ref{42}) gives
\begin{equation}\label{15}
f^L_2(m)=(m-1)f^L_2(2)-(m-2)f^L_2(1)\ \ {\rm if}\ \ m\neq0.
\end{equation}
Taking $m=-1$ in (\ref{15}), we have
\begin{equation}\label{53}
 f^L_2(-1)=-2f^L_2(2)+3f^L_2(1).
\end{equation}
Furthermore, letting $\la=-1$, $m=1$, $n=-1$ in (\ref{8}), we have
$f^L_2(1)+f^L_2(-1)=2f^L_2(0)$. Using this in (\ref{53}), one has $f^L_2(0)=-f^L_2(2)+2f^L_2(1)$.
Thus (\ref{15}) holds for all $m\in\Z$. Taking $m=-2$ in (\ref{15}), we obtain
\begin{equation}\label{004}
f^L_2(-2)=-3f^L_2(2)+4f^L_2(1).
\end{equation}
Taking $\la=-1$ in (\ref{main}), then using (\ref{53}) and (\ref{004}), one can deduce $f^L_2(2)=2f^L_2(1)$. Thus (\ref{15}) gives  $f^L_2(m)=mf^L_2(1)$
for all $m\in\Z$. Taking $\la=-1$ in (\ref{08}) and using $f^L_2(-1)=-f^L_2(1)$, one can deduce $f^L_2(1)=0$. Thus $f^L_2(n)=0$ for all $n\in\Z$.

\begin{subcase}
$\la=1$.
\end{subcase}
Letting $\la=1$ in (\ref{8}), we obtain $(m-n)\big(f^L_2(m+n)-f^L_2(n)-f^L_2(m)\big)=0$.
Utilizing the similar technique to that of (1), one can deduce $f^L_2(m)=mf^L_2(1)$ for all $m\in\Z$.

\begin{subcase}
$\la\notin\{-5,-3,\pm1\}$.
\end{subcase}
Taking $n=0$ in (\ref{8}), we have $(\la+1)mf^L_2(0)=0$ for all $m\in\Z$. This forces $f^L_2(0)=0$. Using this and replacing $n$ by $-m$ in (\ref{8}), we obtain $(\la+3)\big(f^L_2(m)+f^L_2(-m)\big)=0$. Thus $f^L_2(m)=-f^L_2(-m)$ for all $m\in\Z$ since $\la\neq-3$. Thus
\begin{equation}\label{22}
f^L_2(1)=-f^L_2(-1),\ \ f^L_2(2)=-f^L_2(-2).
\end{equation}
Taking $n=2$ in (\ref{8}), we have
\begin{eqnarray}\label{9}
2(m-\la-1)f^L_2(m)-\big(4-(\la+1)m\big)f^L_2(2)=2(m-2)f^L_2(m+2).
\end{eqnarray}
Replacing $m$ by $m+2$, $n$ by $-2$ in (\ref{8}) and using (\ref{22}), we have
\begin{equation}\label{10}
2(m+\la+3)f^L_2(m+2)-\big(4+(\la+1)(m+2)\big)f^L_2(2)=2(m+4)f^L_2(m).
\end{equation}
Using (\ref{9}) and (\ref{10}), one can deduce $(\la+5)(\la-1)\big(2f^L_2(m)-mf^L_2(2)\big)=0$, which gives
$f^L_2(m)=\frac{1}{2}f^L_2(2)m$. Taking $m=2$, $n=-1$ in (\ref{8}) and using (\ref{22}), we can deduce 
$f^L_2(2)=2f^L_2(1)$. Thus $f^L_2(m)=mf^L_2(1)$ for all $m\in\Z$. We have completed the proof of (2).

Next we begin the proof of (3) and (4) of this lemma. Applying $d$ to $[L_n,L_m]=(m-n)L_{m+n}$, then comparing the coefficients of $M_{m+n}$, one has
\begin{eqnarray}
\label{3}(m-\la n+2\mu)f^L_3(m)-(n-\la m+2\mu)f^L_3(n)=(m-n)f^L_3(m+n).
\end{eqnarray}
\setcounter{case}{0}

(3)\  \ If $\mu\notin\Z$, then taking $n=0$ in (\ref{3}), we have
$2\mu f^L_3(m)=(2\mu-\la m)f^L_3(0)$, which gives
$$
f^L_3(m)=\frac{1}{2\mu}(2\mu-\la m)f^L_3(0).
$$

(4)\  \ Taking $\mu=0$ in (\ref{3}), we have
\begin{eqnarray}
&&\label{28}(m-\la n)f^L_3(m)-(n-\la m)f^L_3(n)=(m-n)f^L_3(m+n).
\end{eqnarray}
Taking $n=1,2$ in (\ref{28}), we have
\begin{eqnarray}
\label{30}&& (m-\la)f^L_3(m)-(1-\la m)f^L_3(1)=(m-1)f^L_3(m+1),\\
\label{31}&& (m-2\la)f^L_3(m)-(2-\la m)f^L_3(2)=(m-2)f^L_3(m+2).
\end{eqnarray}
Setting $n=1$ and replacing $m$ by $m+1$ in (\ref{28}), one can deduce
\begin{eqnarray}\label{32}
 (m+1-\la)f^L_3(m+1)-\big(1-\la(m+1)\big)f^L_3(1)=mf^L_3(m+2).
\end{eqnarray}
Thus the three equations (\ref{30})--(\ref{32}) imply
\begin{eqnarray}
&&\nonumber\ \ \ (\la-1) \big( (\la+2) m-2\la\big)f^L_3(m)\\
&&\label{33}=m (m-2)(2 - \la + \la^2 - 2 \la m)f^L_3(1)-m(m-1)(2-\la m)f^L_3(2).
\end{eqnarray}
\setcounter{subcase}{0}
\begin{case}
$\la= -2$.
\end{case}
Setting $\la=-2$ in (\ref{33}), one has
\begin{eqnarray*}\label{39}
f^L_3(m)=\frac{1}{6}m(m+1)(m-1)f^L_3(2)-\frac{1}{3}m (m-2)(m+2)f^L_3(1),\ \ \forall\ m\in\Z.
\end{eqnarray*}

\begin{case}
$\la= -1$.
\end{case}
In this case (\ref{33}) gives
\begin{eqnarray}\label{46}
f^L_3(m)=\frac{1}{2}m(m-1)f^L_3(2)-m (m-2)f^L_3(1)\ \ {\rm for}\ \ m\neq-2.
\end{eqnarray}
Taking $\la=-1$, $m=-2$ and $n=1$ in (\ref{28}), one has
\begin{equation}\label{47}
f^L_3(-2)=3f^L_3(-1)+f^L_3(1).
\end{equation}
Setting $m=-1$ in (\ref{46}), we have $f^L_3(-1)=f^L_3(2)-3f^L_3(1)$. Using this in (\ref{47}), we obtain $f^L_3(-2)=3f^L_3(2)-8f^L_3(1)$. Thus (\ref{46}) holds for all $m\in\Z$.
\begin{case}
$\la=1$.
\end{case}
In this case (\ref{28}) gives $(m-n)\big(f^L_3(m+n)-f^L_3(m)-f^L_3(n)\big)=0$. Using the similar discussions to the proof of (1), we have $f^L_3(m)=mf^L_3(1)$ for all $m\in\Z$.
\begin{case}
$\la\notin\{-2,\pm1,0\}$.
\end{case}
Letting $n=1$ in (\ref{28}), we have
\begin{eqnarray}
\label{34}&& (m-\la)f^L_3(m)-(1-\la m)f^L_3(1)=(m-1)f^L_3(m+1).
\end{eqnarray}
Replacing $m$ by $m+1$, $n$ by $-1$ in (\ref{28}), we obtain
\begin{equation}\label{35}
 (m+1+\la)f^L_3(m+1)+\big(1+\la( m+1)\big)f^L_3(-1)=(m+2)f^L_3(m).
\end{equation}
Thus, using $(\ref{34})$ and $(\ref{35})$, one can deduce
\begin{equation}\label{36}
(1 -\la) (2 + \la)f^L_3(m)=(1-\la m)(m+1+\la)f^L_3(1)-(m-1)\big(\la(m+1)+1\big)f^L_3(-1).
\end{equation}
Taking $m=-1$ in (\ref{36}), we have
$\la(1+\la)\big(f^L_3(1)+f^L_3(-1)\big)=0$, which forces $f^L_3(1)=-f^L_3(-1)$. Then (\ref{36}) gives
$( 1-\la ) ( \la+2) \big(f^L_3(m)-mf^L_3(1)\big)=0$. Thus $f^L_3(m)=mf^L_3(1)$ for all $m\in\Z$.

Hence denoted by $f^L_1(1)=a$, $f^L_2(1)=b$, $f^L_2(0)=e$,
$f^L_3(0)=\bar{e}$, $f^L_3(2)=c$ and $f^L_3(1)=\bar{c}$, the lemma follows.\QED

\begin{lemm}\label{lemma3}  Define $a$, $b$ and $e$ as those given in Lemma \ref{lemma2}. For any $n\in\Z$ and some $\bar{a},\,\hat{b},\,{\bar{b}}\in\C$, we have

(1)\ \ $f_1^Y(n)=0$.

(2)\ \ $f_2^Y(n)=an+\bar{a}$.

(3)\ \ If $\mu\notin\Z$, then $f_3^Y(n)=-\frac{e}{\mu}n$.

(4)\ \ If $\mu=0$, then

\ \ \ \ \ \ \ \ \ \ \ \ \ \ \ \ \ \ $f_3^Y(n)=\left\{\begin{array}{llllllll}%
\hat{b}n&{\rm if}\ \ \la=-1,\vs{6pt}\\
bn+{\bar{b}}&{\rm if}\ \ \la=1,\vs{6pt}\\
\frac{2b}{1+\la}n&{\rm if}\ \ \la\neq\pm1.
\end{array}\right.$
\end{lemm}
\ni{\bf Proof}\ \ (1)\ \ Applying $d$ to $[L_n,Y_m]=(m-\frac{1+\la}{2}n+\mu)Y_{m+n}$ and comparing the coefficients of $L_{m+n}$, we have
\begin{eqnarray}
\label{600}2(m-n)f^Y_1(m)-\big(2m-(1+\la)n+2\mu\big)f^Y_1(n+m)=0.
\end{eqnarray}
\setcounter{case}{0}
\begin{case}
$\mu\notin\Z$.
\end{case}
Setting $n=0$ in (\ref{600}), we have $\mu f^Y_1(m)=0$, which gives $f^Y_1(m)=0$ for all $m\in\Z$.
\begin{case}
$\mu=0$.
\end{case}
In this case (\ref{600}) gives
\begin{eqnarray}
\label{60}2(m-n)f^Y_1(m)-\big(2m-(1+\la)n\big)f^Y_1(n+m)=0.
\end{eqnarray}
Setting $m=0$ in (\ref{60}), we have
\begin{eqnarray}\label{63}
(1+\la)nf^Y_1(n)=2nf^Y_1(0).
\end{eqnarray}
\setcounter{subcase}{0}
\begin{subcase}
$\la\neq\pm1$.
\end{subcase}
By (\ref{63}), we have
\begin{equation}\label{603}
nf^Y_1(n)=\frac{2}{1+\la}nf^Y_1(0),\ \ \forall\ n\in\Z.
\end{equation}
Multiplying (\ref{60}) by $m(n+m)$ and using (\ref{603}),  we deduce
$(\la-1) m n (m + n)f^Y_1(0)=0$ for all $m$, $n\in\Z$, since $\la\neq1$, which forces $f^Y_1(0)=0$. Using this in (\ref{603}), we obtain $f^Y_1(n)=0$ for all $n\in\Z$.
\begin{subcase}
$\la=-1$.
\end{subcase}
Setting $\la=-1$ in (\ref{60}), we have
\begin{equation}\label{66}
(m-n)f^Y_1(m)-mf^Y_1(n+m)=0.
\end{equation}
Replacing $n$ by $-m$ in (\ref{66}), we have $f^Y_1(m)=\frac{1}{2}f^Y_1(0)$ for $m\neq0$. Furthermore, letting $\la=-1$ in (\ref{63}), we have $nf^Y_1(0)=0$ for all $n\in\Z$, which forces $f^Y_1(0)=0$. Thus $f^Y_1(m)=0$ for all $m\in\Z$.
\begin{subcase}
$\la=1$.
\end{subcase}
Letting $\la=1$ in (\ref{60}), we have
\begin{eqnarray}
\label{601}(m-n)\big(f^Y_1(m)-f^Y_1(n+m)\big)=0.
\end{eqnarray}
Setting $m=0$ in (\ref{601}), we obtain $n\big(f^Y_1(n)-f^Y_1(0)\big)=0$, which gives
\begin{equation}\label{605}
f^Y_1(n)=f^Y_1(0),\ \ \forall\ n\in\Z.
\end{equation}

Applying $d$ to $[Y_2,M_3]=0$ and comparing the coefficients of $M_{5}$, we have
\begin{equation}\label{602}
f^M_2(3)=-f^Y_1(2)=-f^Y_1(0).
\end{equation}
Furthermore,  applying $d$ to $[Y_1,Y_2]=M_3$ and comparing the coefficients, we have $f^M_2(3)=f^Y_1(1)+f^Y_1(2)$. This together with (\ref{605}) and (\ref{602}) gives $f^Y_1(0)=0$.  Thus $f^Y_1(n)=0$ for all  $n\in\Z$.

(2)\ \ From Lemma {\ref{lemma2}}, we have $f^L_1(n)=an$ for some $a\in\C$. Thus applying $d$ to $[L_n,Y_m]=(m-\frac{1+\la}{2}n+\mu)Y_{m+n}$ and comparing the coefficient of $Y_{m+n}$, we have
\begin{eqnarray}
\label{61}(2m-(1+\la)n+2\mu)\big(f^Y_2(m)-f^Y_2(n+m)+an\big)=0.
\end{eqnarray}
\setcounter{case}{0}
\begin{case}
$\la=-1$.
\end{case}
In this case (\ref{61}) gives
\begin{equation}\label{68}
(m+\mu)\big(f^Y_2(m)-f^Y_2(n+m)+an\big)=0.
\end{equation}
Recalling that $\mu$ satisfies (\ref{mu}), taking $m=1$ and replacing $n$ by $n-1$ in (\ref{68}), we obtain
\begin{equation}\label{75}
f^Y_2(n)=a(n-1)+f^Y_2(1),\ \ \forall\,n\in\Z.
\end{equation}
Taking $n=0$ in (\ref{75}), we obtain $f^Y_2(1)=f^Y_2(0)+a$. Thus (\ref{75}) gives $f^Y_2(n)=an+f^Y_2(0)$ for all $n\in\Z$.
\begin{case}
$\la\neq-1$.
\end{case}
If $\mu=0$, then setting $m=0$ in (\ref{61}), we obtain
$(1+\la)n\big(f^Y_2(0)-f^Y_2(n)+an\big)=0$,
which gives $f^Y_2(n)=an+f^Y_2(0)$ for all $n\in\Z$.

Suppose $\mu\notin\Z$. Replacing  $n$ by $-m$ in (\ref{61}), we obtain
\begin{equation}\label{69}
\big((3+\la)m+2\mu)\big(f^Y_2(m)-am-f^Y_2(0)\big)=0.
\end{equation}

If $(3+\la)m+2\mu\neq0$ for all $m\in\Z$, then (\ref{69}) gives $f^Y_2(m)=am+f^Y_2(0)$.

If $(3+\la)m'+2\mu=0$ for some $m'\in\Z$, then replacing $m$ by $m'+1$ in (\ref{69}), we have
\begin{equation}\label{add}
(3+\la)\big(f^Y_2(m'+1)-a(m'+1)-f^Y_2(0)\big)=0.
\end{equation}
It is obvious that $\la\neq-3$, otherwise $\mu=0$. Thus (\ref{add}) gives
\begin{equation}\label{0014}
f^Y_2(m'+1)=a(m'+1)+f^Y_2(0).
\end{equation}
Choosing $n=1$ and replacing $m$ by $m'$ in (\ref{61}), we obtain
\begin{equation}\label{0015}
(2\mu-\la+2m'-1)\big(a+f^Y_2(m')-f^Y_2(m'+1)\big)=0.
\end{equation}
Combining (\ref{0014}) and (\ref{0015}), one can deduce
\begin{equation}\label{0016}
(2\mu-\la+2m'-1)\big(f^Y_2(m')-am'-f^Y_2(0)\big)=0.
\end{equation}
Furthermore, using $2\mu=-(3+\la)m'$ in (\ref{0016}), we have
\begin{equation}\label{0017}
(m'+1)(1+\la)\big(f^Y_2(m')-am'-f^Y_2(0)\big)=0.
\end{equation}
In this case $\la\neq-1$, we obtain
\begin{equation}\label{0018}
(m'+1)\big(f^Y_2(m')-am'-f^Y_2(0)\big)=0.
\end{equation}
If $m'\neq-1$, then (\ref{0018}) gives $f^Y_2(m')=am'+f^Y_2(0)$.

Suppose that $m'=-1$. Then $2\mu=3+\la$. Setting $m=0$, $n=-1$ in (\ref{61}), we obtain
\begin{eqnarray}\label{0019}
(\la+2)\big(f^Y_2(0)-f^Y_2(-1)-a\big)=0.
\end{eqnarray}
If $\la\neq-2$, then (\ref{0019}) forces $f^Y_2(-1)=-a+f^Y_2(0)$.
If $\la=-2$, then setting $\la=-2$, $m=-1$ and  $n=1$ in (\ref{61}), we obtain
\begin{eqnarray}\label{0020}
{(2\mu-1)\big(a+f^Y_2(-1)-f^Y_2(0)\big)=0,}
\end{eqnarray}
since $\mu$ satisfies (\ref{mu}), which forces $f^Y_2(-1)=-a+f^Y_2(0)$.
Hence $f^Y_2(m)=am+f^Y_2(0)$ for all $m\in\Z$. By now we have completed the proof of (2).

Next we begin the proof of (3) and (4) of this lemma. Applying $d$ to $[L_n,Y_m]=(m-\frac{1+\la}{2}n+\mu)Y_{m+n}$ and comparing the coefficients of $M_{m+n}$, we have
\begin{eqnarray}\label{620}
2(m-n)f^L_2(n)+2(m-\la n+2\mu)f^Y_3(m)=(2m-(1+\la)n+2\mu)f^Y_3(n+m).
\end{eqnarray}
\setcounter{case}{0}

(3)\ \  If $\mu\notin\Z$, then taking $n=0$ in (\ref{620}), we obtain $f^Y_3(m)=-\displaystyle\frac{e}{\mu}m$ for $e=f^L_2(0)$.

(4)\ \  Taking $\mu=0$ in (\ref{620}), we obtain
\begin{eqnarray}\label{62}
2(m-n)f^L_2(n)+2(m-\la n)f^Y_3(m)=(2m-(1+\la)n)f^Y_3(n+m).
\end{eqnarray}
By Lemma {\ref{lemma2}} (2), one has $f^L_2(n)=bn(1-\delta_{\la,-1})$ for some $b\in\C$ .
\setcounter{subcase}{0}
\begin{case}
$\la=-1$.
\end{case}
In this case (\ref{62}) gives
\begin{equation}\label{78}
(m+n)f^Y_3(m)=mf^Y_3(m+n).
\end{equation}
Taking $m=1$ and replacing $n$ by $n-1$ in (\ref{78}), we have $f^Y_3(n)=nf^Y_3(1)$ for all $n\in\Z$.
\begin{case}
$\la=1$.
\end{case}
Letting $\la=1$ in (\ref{62}), we obtain
 \begin{eqnarray}\label{007}
(m-n)\big(f^Y_3(n+m)-f^Y_3(m)-{b}n\big)=0.
\end{eqnarray}
Thus (\ref{007}) gives $f^Y_3(n)={b}n+f^Y_3(0)$ for all $n\in\Z$.
 \begin{case}
$\la\neq\pm1$.
\end{case}
Setting $m=0$ in   (\ref{62}), we have
\begin{equation}\label{74}
f^Y_3(n)=\frac{2}{1+\la}\big(bn+\la f^Y_3(0)\big)\ \ {\rm for}\ \ n\neq0.
\end{equation}
Furthermore, taking $n=\pm1$ in (\ref{74}) respectively, we have
\begin{equation}\label{005}
f^Y_3(1)=\frac{2}{1+\la}\big(\la f^Y_3(0)+b\big),\ \ f^Y_3(-1)=\frac{2}{1+\la}\big(\la f^Y_3(0)-b\big).
\end{equation}
Taking $n=1$, $m=-1$ in (\ref{62}) and using (\ref{005}), one can deduce
$(\la-1)f^Y_3(0)=0$, which forces $f^Y_3(0)=0$ since $\la\neq1$. Then (\ref{74}) implies $f^Y_3(n)=\frac{2b}{1+\la}n$ for all $n\in\Z$.

Hence, denoted by $f^Y_2(0)=\bar{a}$, $f^Y_3(1)=\hat{b}$ and $f^Y_3(0)={\bar{b}}$, the lemma follows.
\QED
\begin{lemm}\label{lemma4}  Let $a$ and $\bar{a}$ be as those given in Lemmas \ref{lemma2} and \ref{lemma3} respectively. For any $n\in\Z$, we have
$$f_1^M(n)=f^M_2(n)=0,\ \ f_3^M(n)=an+2\bar{a}.$$
\end{lemm}
\ni{\bf Proof}\ \ \ It follows from Lemmas \ref{lemma2}  and \ref{lemma3} that
there exists some $a$ and $\bar{a}$ such that
\begin{eqnarray*}
f^Y_1(n)=0,\ \ f^Y_2(n)=an+\bar{a},\ \ \forall\,n\in\Z.
\end{eqnarray*}
Applying $d$ to $[Y_n,Y_m]=(m-n)M_{m+n}$ and comparing the coefficients of $L_{m+n}$, $Y_{m+n}$ and $M_{m+n}$ respectively, we have
\begin{eqnarray}
\label{85}&&(m-n)f_1^M({m+n})=(m-n)f_2^M({m+n})=0,\\
\label{87}&&(m-n)\big(f^M_3(m+n)-a(n+m)-2\bar{a}\big)=0.
\end{eqnarray}
Setting $m=0$ and replacing $m$ by $-n$ in  (\ref{85}), we obtain
$$nf_1^M(n)=nf_2^M(n)=0,\ \ nf_1^M(0)=nf_2^M(0)=0,\ \ \forall\ n\in\Z.$$
This gives $f_1^M(n)=f_2^M(n)=0$ for all $n\in\Z$.

Taking $m=0$ in (\ref{87}), one can deduce $f^M_3(n)=an+2\bar{a}$ for $n\neq0$. Taking $m=1$ and $n=-1$ in (\ref{87}), one can deduce $f^M_3(0)=2\bar{a}$. Thus $f^M_3(n)=an+2\bar{a}$ for all $n\in\Z$.\QED

We construct some possible outer derivations of $\mathscr{L}$. Under the condition $\mu=0$, for each $\la\in\{-2,\pm1\}$, the following maps $D_{\la}$ or $\overline{D}_\la$ defined by
\begin{eqnarray*}
&&D_{-2}(L_n)=n^3 M_n,\ D_{-1}(L_n)=n^2 M_n,\ \overline{D}_{-1}(Y_n)=nM_n,\ D_{1}(Y_n)=M_n,
\end{eqnarray*}
are outer derivations of $\mathscr{L}$, where all other terms are vanishing.
Besides, we can define another outer derivation $D$ of $\mathscr{L}$, which does not depend on $\la$ and $\mu$:
\begin{eqnarray*}
D: L_n\mapsto 0,\  Y_n\mapsto Y_n,\ M_n\mapsto 2M_n.
\end{eqnarray*}
It  is easy to verity that  for each $\la\in\{-2,\pm1\}$, $D$, $D_{\la}$ and $\overline{D}_\la$ are linearly independent.
\begin{theo} \label{theorem1}\ \
(1)\ \ If $\mu\notin\frac{1}{2}\Z$,  or $\mu=0$ but $\la\notin\{-2,0,\pm1\}$, then Der$\mathscr{L}=ad\mathscr{L}\oplus \C D$ .

(2)\ \ If $\mu=0$, then
\begin{eqnarray}\label{Der}
Der\mathscr{L}=\left\{\begin{array}{lllllllll}
 ad\mathscr{L}\oplus\C D\oplus \C D_{-2},&\la=-2,\vs{4pt}\\
 ad\mathscr{L}\oplus\C D\oplus \C D_{-1}\oplus \C \overline{D}_{-1},&\la=-1,\vs{4pt}\\
ad\mathscr{L}\oplus \C D\oplus\C {D}_{1},&\la=1.
\end{array}\right.
\end{eqnarray}
\end{theo}
\ni{\bf Proof}\ \ \ Take $d_0\in (Der\mathscr{L})_0$.

(1)\ \ Suppose $\mu\notin\frac{1}{2}\Z$. It follows from Lemmas \ref{lemma2}--\ref{lemma4} that there exist $a$, $\bar{a}$, $e$ and $\bar{e}\in\C$ such that
\begin{eqnarray*}
&&d_0(L_n)=anL_n+\frac{e}{2\mu}\big(2\mu-(\la+1)n\big)Y_n+\frac{\bar{e}}{2\mu}(2\mu-\la n)M_n,\\
&&d_0(Y_n)=\big(an+\bar{a}\big)Y_n-\frac{e}{\mu}nM_n,\ \ d_0(M_n)=(an+2\bar{a})M_n.
\end{eqnarray*}
Set $\al=aL_0-{e}{\mu^{-1}}Y_{0}-{\bar{e}}{(2\mu)^{-1}}M_0$, then $d_0=ad_{\al}+({\bar{a}-a\mu})D$.

If $\mu=0$ and $\la\notin\{-2,0,\pm1\}$, by  Lemmas \ref{lemma2}--\ref{lemma4}, there exist $a'$, ${\bar{a}'}$, $b'$ and $e'\in\C$ such that
\begin{eqnarray*}
&&d_0(L_n)=a'nL_n+b'nY_n+e'nM_n,\\
&&d_0(Y_n)=(a'n+\bar{a}')Y_n+\frac{2b'}{1+\la}nM_n,\ \ d_0(M_n)=(a'n+2\bar{a}')M_n.
\end{eqnarray*}
Set $\beta=a'L_0+{2b'}{(1+\la)^{-1}}Y_0+{e'}{\la^{-1}}M_0$. Then $d_0=ad_{\beta}+\bar{a}'D$.

(2)\ \ We shall divide the proof of (\ref{Der}) into the following four cases.
\setcounter{case}{0}
\begin{case}
$\la=-2$.
\end{case}
By  Lemmas \ref{lemma2}--\ref{lemma4}, there exist $a_1$, ${\bar{a}_1}$,  $b_1$, $c_1$ and $\bar{c}_1\in\C$ such that
\begin{eqnarray*}
&&d_0(L_n)=a_1nL_n+b_1nY_n+\big(\frac{c_1}{6}(n^3-n)-\frac{\bar{c}_1}{3} (n^3-4n)\big)M_n,\\
&&d_0(Y_n)=(a_1n+\bar{a}_1)Y_n-2b_1nM_{n},\ \ d_0(M_n)=(a_1n+2\bar{a}_1)M_n.
\end{eqnarray*}
Set $\al_1=a_1L_0-2b_1Y_0+\frac{1}{12}(c_1-8\bar{c}_1)M_0$. Then we obtain $d_0=ad_{\al_1}+\bar{a}_1D+\frac{1}{6}(c_1-2\bar{c}_1)D_{-2}$.
\begin{case}
$\la=-1$.
\end{case}
By  Lemmas \ref{lemma2}--\ref{lemma4}, there exist $a_2$, ${\bar{a}_2}$, ${b}_2$, $c_2$ and $\bar{c}_2\in\C$ such that
\begin{eqnarray*}
&&d_0(L_n)=a_2nL_n+\big(\frac{c_2}{2}(n^2-n)-\bar{c}_2(n^2-2n)\big)M_n,\\
&&d_0(Y_n)=(a_2n+\bar{a}_2)Y_n+{b}_2nM_n,\ \ d_0(M_n)=(a_2n+2\bar{a}_2)M_n.
\end{eqnarray*}
Set $\al_2=a_2L_0+\frac{1}{2}(c_2-4\bar{c}_2)M_0$. Then $d_0=ad_{\al_2}+\bar{a}_2D+\frac{1}{2}(c_2-2\bar{c}_2)D_{-1}+{b}_2\overline{D}_{-1}$.
\begin{case}
$\la=1$.
\end{case}
By  Lemmas \ref{lemma2}--\ref{lemma4}, there exist $a_4$, ${\bar{a}_4}$, $b_4$, $\bar{b}$ and $c_4\in\C$ such that
\begin{eqnarray*}
&&d_0(L_n)=a_4nL_n+{b_4}nY_n+c_4nM_n,\\
&&d_0(Y_n)=(a_4n+\bar{a}_4)Y_n+({b_4}n+\bar{b})M_n,\ \ d_0(M_n)=(a_4n+2\bar{a}_4)M_n.
\end{eqnarray*}
Set $\al_4=a_4L_0+{b_4}Y_0+c_4M_0$. Then $d_0=ad_{\al_4}+\bar{a}_4D+\bar{b}{D}_1$.

Hence, the theorem follows from Lemma \ref{lemma1}.\QED

\section{Automorphism group of $\mathscr{L}$}
\setcounter{section}{3}\setcounter{equation}{0}

\indent\ \ \ \ \ In this section we denote by $\C^*$ the set of all nonzero complex numbers and $Aut\mathscr{L}$ and $Inn\mathscr{L}$ the sets of automorphisms and inner automorphisms of $\mathscr{L}$. Firstly we need to introduce some technical lemmas.

\begin{lemm}\label{lemma5}
Let  $\sigma\in Aut\mathscr{L}$, $\epsilon\in\{\pm1\}$ with $\epsilon=1$ for $\mu\notin\frac{1}{2}\Z$. Then there exist some $\al$, $\beta\in\C^*$, $y_n\in \mathcal{Y}$ and $m_n$, $m'_n\in \mathcal{M}$ such that

(1)\ \ $\sigma(L_n)=\epsilon\al^n L_{\epsilon n}+y_n+m_n$,

(2)\ \ $\sigma(Y_n)=\al^n\beta Y_{\epsilon n}+m'_n$,

(3)\ \ $\sigma(M_n)=\epsilon\al^n\beta^2M_{\epsilon n}$.

\end{lemm}
\ni{\bf Proof}\ \ (1)\ \ Note that $\sigma|_{\mathcal{L}}$ is an automorphism of the  Witt algebra, so there exist some $\epsilon\in\{\pm1\}$ and $\al\in\C^*$ such that
\begin{eqnarray}\label{401}
\sigma(L_n)=\epsilon\al^n L_{\epsilon n}+y_n+m_n\ \ \ {\rm for\ some}\  y_n\in \mathcal{Y}, \ m_n\in \mathcal{M}.
\end{eqnarray}

(2)\ \ Since $\mathcal{I}=\mathcal{Y}\oplus \mathcal{M}$ is the unique maximal idea of $\mathscr{L}$, one can write
\begin{eqnarray}\label{402}
\sigma(Y_n)=\mbox{$\sum\limits_{i\in S}^{}$}b_{n_i}Y_{i}+m'_n,
\end{eqnarray}
for some $b_{n_i}\in \C^*$, $m'_n\in \mathcal{M}$ and $S\subseteq\Z$. By (\ref{401}), there exist some $y_0\in \mathcal{Y}$, $m_0\in \mathcal{M}$ such that $\sigma(L_0)=\epsilon L_{0}+y_0+m_0$. Then applying $\sigma$ to $(n+\mu)Y_{n}=[L_0,Y_n]$, we have
\begin{eqnarray*}\label{40201}
&&(n+\mu)\Big(\mbox{$\sum\limits_{i\in S}^{}$}b_{n_i}Y_{i}+m'_n\Big)=\big[\epsilon L_{0}+y_0+m_0,\mbox{$\sum\limits_{i\in S }^{}$}b_{n_i}Y_{i}+m'_n]=\mbox{$\sum\limits_{i\in S }^{}$}\epsilon b_{n_i}(i+\mu)Y_{i}+m''_n,
\end{eqnarray*}
for some $m''_n\in \mathcal{M}$. Comparing the coefficients of $Y_i$, we obtain
\begin{equation}\label{0023}
b_{n_i}\big(i-\epsilon n-(\epsilon-1)\mu\big)=0,\ \forall\,\,i\in S.
\end{equation}
If $\mu=0$, then (\ref{0023}) gives $i=\epsilon n$ for all $i\in S$. If $\mu\notin\frac{1}{2}\Z$, then (\ref{0023}) forces $\epsilon=1$, in which case (\ref{0023}) gives $i=n$ for all $i\in S$. Hence, by (\ref{402}), we can write
\begin{equation}\label{y}
\sigma(Y_n)=b_nY_{\epsilon n}+m'_n,
\end{equation}
for some $b_n\in\C^*$ with $\epsilon=1$ for $\mu\notin\frac{1}{2}\Z$.

Applying $\sigma$ to $(m-\frac{\la+1}{2}n+\mu)Y_{m+n}=[L_n,Y_m]$ and using (\ref{401}) and (\ref{y}), we obtain
\begin{eqnarray*}
&&\!\!\!\!\!\!(m-\frac{\la+1}{2}n+\mu)\big(b_{m+n}Y_{\epsilon(m+n)}+m'_{m+n}\big)\\
&\!\!\!\!\!\!\!\!\!\!\!\!\!\!\!\!\!\!\!=\!\!\!\!\!\!\!\!\!&\!\!\!\!\!\![\epsilon\al^n L_{\epsilon n}+y_n+m_n,b_mY_{\epsilon m}+m'_m]=\al^n b_m(m-\frac{\la+1}{2}n+\mu)Y_{\epsilon (m+n)}+\bar{m}_n,
\end{eqnarray*}
for some $\bar{m}_n\in \mathcal{M}$. Thus we have
\begin{eqnarray*}
(m-\frac{\la+1}{2}n+\mu)\big((b_{m+n}-\al^n b_m)Y_{\epsilon (m+n)}+m'_{m+n}\big)-\bar{m}_n=0,
\end{eqnarray*}
which together with the fact $m'_{m+n}$, $\bar{m}_n\in \mathcal{M}$, gives
\begin{eqnarray}\label{0039}
\big(2m-(\la+1)n+2\mu\big)(b_{m+n}-\al^nb_m)=0.
\end{eqnarray}
Taking $m=0$ in (\ref{0039}), one has
\begin{eqnarray}
\label{0040}\big(2\mu -(\la+1)n\big)(b_{n}-\al^nb_0)=0.
\end{eqnarray}
If $2\mu-(\la+1)n'=0$ for some $n'\in\Z$, then replacing $n$ by $-n'$, $m$ by $n'$ in (\ref{0039}), we have
 \begin{eqnarray}
\label{0041}( \la +2)n'(b_{n'}-\al^{n'}b_0)=0.
\end{eqnarray}

If $\la\neq-2$, then (\ref{0041}) gives
$b_{n'}=\al^{n'}b_0$.

If $\la=-2$, then $2\mu+n'=0$, since $\mu$ satisfies (\ref{mu}), which forces $n'=0$. Thus
\begin{equation}\label{00412}
b_{n'}=\al^{n'}b_0\ \ \ \ \ {\rm if}\ \ 2\mu-(\la+1)n'=0.
\end{equation}
Hence, (\ref{0040}) together with (\ref{00412}) gives $b_{n}=\al^{n}b_0$ for all $n\in\Z$. Using this in (\ref{y}), one can write
\begin{equation}\label{yy}
\sigma(Y_n)=\al^{n}b_0Y_{\epsilon n}+m'_n,
\end{equation}
for some $m'_n\in \mathcal{M}$.

(3)\ \ Applying $\sigma$ to $(n-2m)M_{n}=[Y_m,Y_{n-m}]$ and using (\ref{yy}), we have
$(n-2m)\sigma(M_{n})=[\al^m b_0Y_{\epsilon m}+m'_m,\al^{n-m}b_0Y_{\epsilon(n-m)}+m'_{n-m}]$ for some $m'_m,m'_{n-m}\in \mathcal{M}$, which gives
\begin{equation}\label{new}
(n-2m)\big(\sigma(M_{n})-\epsilon\al ^{n}b_0^2 M_{\epsilon n}\big)=0.
\end{equation}
Letting $m=0$ in (\ref{new}), we have $\sigma(M_{n})=\epsilon\al ^{n}b_0^2 M_{\epsilon n}$ for $n\neq0$. Furthermore, letting $n=0$ in (\ref{new}), one has $\sigma(M_{0})=\epsilon b_0^2 M_{0}$. Thus we obtain $\sigma(M_{n})=\epsilon\al ^{n}b_0^2M_{\epsilon n}$.

 Hence, denoted by $b_0=\beta$, the lemma follows.\QED

\begin{lemm}\label{lemma9}
Let $\bar{f}$, $f$ and $g$ be $\C$-linear maps from $\Z$ to $\C$. Define the $\C$-linear map $\phi:\mathscr{L}\rightarrow\mathscr{L}$ by $$\phi(L_n)=L_n+\bar{f}(n)Y_n+f(n)M_n,\ \phi(Y_n)=Y_n+g(n) M_n,\ \phi(M_n)=M_n.$$
If $\mu=0$ and $\phi\in Aut\mathscr{L}$, then there exist some $a$, $\bar{a}$, $b$, $c$, $\bar{c}\in\C$ such that

(1)\ \ $$\bar{f}(n)=bn(1-\delta_{\la,-1}),\ \ \ g(n)=\left\{\begin{array}{llllllll}
an&{\rm if}\ \ \la=-1,\vs{6pt}\\
bn+\bar{a}&{\rm if}\ \ \la=1,\vs{6pt}\\
\displaystyle\frac{2b}{1+\la}n&{\rm if}\ \ \la\neq\pm1;
\end{array}\right.$$

(2)\ \ \begin{eqnarray*}
&&\label{f}f(n)=\left\{\begin{array}{llllllll}%
\displaystyle\frac{1}{6}c(n^3-n)-\frac{1}{3} \bar{c}(n^3-4n)+\displaystyle\frac{1}{3}b^2n(n-2)(n-1)&{\rm if}\ \ \la=-2,\vs{6pt}\\
\displaystyle\frac{1}{2}c(n^2-n)-\bar{c}(n^2-2n)&{\rm if}\ \ \la=-1,\vs{6pt}\\
\bar{c}n+\displaystyle\frac{b^2}{1+\la}n(n-1)&{\rm if}\ \ \la\notin\{-2,0,-1\}.
\end{array}\right.
\end{eqnarray*}
\end{lemm}
{\ni\bf Proof}\ \ Applying $\phi$ to $[L_n,L_m]=(m-n)L_{n+m}$ and comparing the coefficients of $Y_{n+m}$,  $M_{n+m}$, we obtain
\begin{eqnarray}
&&\label{0054}\big(2m-(1+\la)n\big)\bar{f}(m)-\big(2n-(1+\la)m\big)\bar{f}(n)
=2(m-n)\bar{f}({m+n}),\\
&&\label{0055}(m-\la n)f(m)-(n-\la m)f(n)=(m-n)\big(f({m+n})-\bar{f}(n)\bar{f}(m)\big).
\end{eqnarray}
Applying $\phi$ to $[L_n,Y_m]=(m-\frac{\la+1}{2}n)Y_{n+m}$ and comparing the coefficients of $M_{n+m}$, we obtain
\begin{eqnarray}
&&\label{0068}2(m-n)\bar{f}(n)+2(m-\la n)g(m)=\big(2m-(\la+1)n)g({n+m}).
\end{eqnarray}
If we replace $\bar{f}$ by $f^L_2$, $g$ by $f^Y_3$ in (\ref{0054}) and (\ref{0068}), then (\ref{0054}) and (\ref{0068}) become (\ref{8}) and (\ref{62}), respectively. Thus if we take $\bar{f}(1)=b$, $g(1)=a$ and  $g(0)=\bar{a}$, then by Lemmas \ref{lemma2} and \ref{lemma3}, (1) follows.

(2)\ \ If $\la=-1$, then noticing  $\bar{f}(n)=0$ in (\ref{0055}), we obtain
\begin{eqnarray}\label{5501}
(m+n)f(m)-(n+m)f(n)=(m-n)f({m+n}).
\end{eqnarray}
If we replace $f^L_3$ by $f$, $\la$ by $-1$  in (\ref{28}), then (\ref{28}) becomes (\ref{5501}). Thus by Lemma \ref{lemma2} (3), we have $f(n)=\frac{1}{2}(n^2-n)f(2)-(n^2-2n)f(1)$.

If $\la\neq-1$, noticing $\bar{f}(n)=bn$ in (\ref{0055}), we have
\begin{eqnarray}\label{5505}
(m-\la n)f(m)-(n-\la m)f(n)=(m-n)\big(f({m+n})-b^2nm\big).
\end{eqnarray}
Taking $n=1$, $2$ in (\ref{5505}) respectively, we have
\begin{eqnarray}
\label{5508}&& (m-\la)f(m)-(1-\la m)f(1)=(m-1)\big(f(m+1)-b^2m\big),\\
\label{5509}&& (m-2\la)f(m)-(2-\la m)f(2)=(m-2)\big(f(m+2)-2b^2m\big).
\end{eqnarray}
Setting $n=1$ and replacing $m$ by $m+1$ in (\ref{5505}), one has
\begin{eqnarray}\label{5510}
(m+1-\la)f(m+1)-\big(1-\la( m+1)\big)f(1)=m\big(f({m+2})-b^2(m+1)\big).
\end{eqnarray}
Thus using (\ref{5508})--(\ref{5510}), one can deduce
\begin{eqnarray}
&&(\la-1) \big(( \la +2)m-2\la\big)f(m)-b^2(\la-2)m(m-2)(m-1)\nonumber\\
&&=m(m-1)(\la m-2)f(2)+ m (m-2)(2 - \la + \la^2 - 2 \la m)f(1)\label{5511}.
\end{eqnarray}
\setcounter{case}{0}
\begin{case}
$\la= -2$.
\end{case}
In this case (\ref{5511}) gives
\begin{eqnarray*}\label{5512}
f(m)=\frac{1}{6}(m^3-m)f(2)-\frac{1}{3}(m^3-4m)f(1)+\frac{b^2}{3}m(m-2)(m-1),\ \ \forall\ m\in\Z.
\end{eqnarray*}
\begin{case}
$\la=1$.
\end{case}
In this case (\ref{5505}) gives
\begin{eqnarray}
&&\label{5502}(m-n)\big(f(m)+f(n)-f({m+n})+b^2nm\big)=0.
\end{eqnarray}
Obviously, $f(0)=0$. Replacing $m$ by $-n$ in (\ref{5502}), one can deduce
\begin{eqnarray}
&&\label{5503}f(-n)+f(n)=b^2n^2,\ \forall\ n\in\Z.
\end{eqnarray}
By (\ref{5502}), we have
\begin{eqnarray}
&&\label{5504}f({m+n})=f(m)+f(n)+b^2nm\ \ {\rm for}\ \ m\neq n.
\end{eqnarray}
Combining (\ref{5503}), (\ref{5504}), we obtain
\begin{eqnarray*}
f(2m)&=&f(3m)+f(-m)-3b^2m^2\\
&=&f(2m)+f(m)+2b^2m^2+f(-2m)+f(m)-2b^2m^2-3b^2m^2\\
&=&2f(m)+b^2m^2.
\end{eqnarray*}
Thus (\ref{5504}) holds for all $m$, $n\in\Z$. Taking $m=1$ and replacing $n$ by $n-1$ in (\ref{5504}), we obtain
$f(n)=f(n-1)+f(1)+b^2(n-1)$. Hence, using induction on $n$, one can deduce $f(n)=nf(1)+\frac{1}{2}b^2n(n-1)$.
\begin{case}
$\la\notin\{-2,0,1\}$.
\end{case}
Taking $m=0$ in (\ref{5505}), one has $\la nf(0)=0$ for all $n\in\Z$. Thus $f(0)=0$. Using this, then setting $m=-1$ and $n=1$ in (\ref{5505}), one can deduce
\begin{eqnarray}\label{5513}
(1+\la )f(-1)=2b^2-(1+\la)f(1).
\end{eqnarray}
Then taking $n=1$ in (\ref{5505}), we have
\begin{eqnarray}
\label{5514}&& (m-1)f(m+1)=(m-\la)f(m)-(1-\la m)f(1)+(m-1)b^2m.
\end{eqnarray}
Setting $n=-1$ and replacing $m$ by $m+1$ in (\ref{5505}), we obtain
\begin{eqnarray}\label{5516}
(m+1+\la)f(m+1)+(1+\la( m+1))f(-1)=(m+2)\big(f({m})+b^2(m+1)\big).
\end{eqnarray}
Multiplying (\ref{5516}) by $(1+\la)(m-1)$, then using (\ref{5513})  and (\ref{5514}), one can deduce
\begin{eqnarray}
&&\nonumber\ \ \ (1-\la)(1+\la)(2+\la)f(m)\\
&&\nonumber=(1-\la)(1+\la)(2+\la)mf(1)+b^2(1-\la) (2+\la)m(m-1),
\end{eqnarray}
since $\la\notin\{-2,\pm1\}$, which gives $f(m)=mf(1)+\frac{b^2}{\la+1}m(m-1)$ for all $m\in\Z$.

Thus denoted by $f(2)=c$ and $f(1)=\bar{c}$, the lemma follows.\QED

\begin{lemm}\label{auto}
(i)\ \  Let $\epsilon\in\{\pm1\}$. If $\mu=0$, then the map
\begin{equation*}
\varphi_\epsilon:\ L_n\mapsto\epsilon L_{\epsilon n},\ Y_n\mapsto Y_{\epsilon n},\  M_n\mapsto \epsilon M_{\epsilon  n},
\end{equation*}
is an automorphism of $\mathscr{L}$. The set $\big\{\varphi_\epsilon\,|\,\epsilon\in\{\pm1\}\big\}\cong\Z_2$ forms a subgroup of $Aut\mathscr{L}$, where $\varphi_\epsilon\varphi_{\epsilon'}=\varphi_{\epsilon\epsilon'}$ for $\epsilon$, $\epsilon'\in\{\pm1\}$.

(ii)\ \ For any $\al$, $\beta\in\C^*$, the map
\begin{equation*}
{\varphi}_{\al,\beta}:\ L_n\mapsto  \al^nL_{n},\ Y_n\mapsto \al^n\beta Y_n,\ M_n\mapsto \al^n\beta^2M_n
\end{equation*}
is an automorphism of $\mathscr{L}$. The set $\{\varphi_{\al,\beta}\,|\,\al,\beta\in\C^*\}\cong\C^*\times\C^*$ forms a subgroup of $Aut\mathscr{L}$, where ${\varphi}_{\al,\beta} {\varphi}_{\al',\beta'}\!=\! {\varphi}_{\al\al',\beta\beta'}$ for $\al$, $\al'$, $\beta$, $\beta'\in\C^*$.

(iii)\ \ For any $b\in\C$, if $\mu=0$, then the map
\begin{eqnarray*}
\phi_{b}(L_n)=\left\{\begin{array}{ll}
L_{n}+bn^3M_n,&\la=-2,\vs{6pt}\\
L_{n}+bn^2M_n,&\la=-1,
\end{array}\right.\phi_{b}(X_n)=X_n\ \ {\rm for}\ X\in\{Y,M\},
\end{eqnarray*}
is an automorphism of $\mathscr{L}$.  The set $\{\phi_{b}\,|\,b\in\C\}\cong\C$ forms a subgroup of $Aut\mathscr{L}$, where $\phi_{b}\phi_{b'}=\phi_{b+b'}$ for $b$, $b'\in\C$.

(iv)\ \ For any $e\in\C$, if $\mu=0$, then the map $\psi_e$ defined by
\begin{eqnarray*}
\begin{array}{rllllll}
&L_n\mapsto  L_{n},&Y_n\mapsto Y_n+enM_n,&M_n\mapsto M_n,&\la=-1,\vs{6pt}\\
&L_n\mapsto  L_{n},&Y_n\mapsto Y_n+eM_n,&M_n\mapsto M_n,&\la=1,
\end{array}
\end{eqnarray*}
 is an automorphism of $\mathscr{L}$. The set $\{\psi_{e}\,|\,e\in\C\}\cong\C$ forms a subgroup of $Aut\mathscr{L}$, where $\psi_{e}\psi_{e'}=\psi_{e+e'}$ for $e$, $e'\in\C$.
\end{lemm}
\ni{\bf Proof}\ \ \ This follows from straightforward verifications, we omit the details here.\QED\vs{6pt}

Introduce the following notation
\begin{equation}\label{inn}
Inn\mathscr{L}={\rm Span}\big\{{\rm exp}(a{\rm ad}L_0+\mbox{$\sum\limits$} b_i{\rm ad}Y_i+\mbox{$\sum\limits$} c_j{\rm ad}M_j)\,|\,a,b_i,c_j\in\C,i,j\in\Z\big\}
\end{equation}
with $a=0$\ \,if\ \,$\mu\notin\Z$.

\begin{theo}\label{theo}
(1)\ \ If $\mu\notin\frac{1}{2}\Z$, then $Aut\mathscr{L}\cong Inn\mathscr{L}\times\C^*\times\C^*$.

(2)\ \ If $\mu=0$, then
\begin{eqnarray*}
Aut\mathscr{L}\cong\left\{\begin{array}{llllll}
Inn\mathscr{L}\times \C^*\times\Z_2\times\C&{\rm if}\ \ \la=-2\ \,{\rm or}\ \,1,\vs{6pt}\\
Inn\mathscr{L}\times \C^*\times\Z_2\times\C\times\C&{\rm if}\ \ \la=-1,\vs{6pt}\\
Inn\mathscr{L}\times \C^*\times\Z_2&{\rm if}\ \ \la\notin\{-2,0,\pm1\}.
\end{array}\right.
\end{eqnarray*}
\end{theo}
\ni{\bf Proof}\ \ \ Let $\sigma$ be an automorphism of $\mathscr{L}$. By Lemma \ref{lemma5}, one can write
\begin{equation}\label{0048}
\sigma(L_0)=\epsilon L_0+b_0Y_0+c_0M_0+\mbox{$\sum\limits_{0\neq i\in S}$}b_iY_i+\mbox{$\sum\limits_{0\neq j\in S'}$}c_jM_j,
\end{equation}
for some $S$, $S'\subset\Z$. Construct an inner automorphism $\theta$ of $\mathscr{L}$
\begin{equation}\label{0049}
\theta={\rm exp}\bigg(\mbox{$\sum\limits_{0\neq i\in S}$}\frac{ib_0b_i}{ (i+\mu)(i+2\mu)}{\rm ad}M_i-\frac{b_i}{\epsilon(i+\mu)}{\,\rm ad\,}Y_i\bigg){\rm exp}\bigg(\mbox{$\sum\limits_{0\neq j\in S'}$}\frac{-c_j}{\epsilon(j+2\mu)}{\,\rm ad}M_j\bigg).
\end{equation}
One can check that   $\theta^{-1}\sigma(L_0)=\epsilon L_0+b_0Y_0+c_0M_0$. Furthermore, since $L_0$ is a semisimple element of $\mathscr{L}$, then $\theta^{-1}\sigma(L_0)$ is also semisimple. If we denote $A$ the matrix of ${\rm ad}\big(\theta^{-1}\sigma(L_0)\big)$ with respect to the basis $\{L_n,Y_n,M_n\}$, then
$$A=\left(
\begin{array}{ccc}
\epsilon n & b_0(\frac{\la+1}{2}n-\mu) & c_0(\la n-2\mu)\\
    0 & \epsilon(n+\mu) & b_0n \\
    0 & 0 & \epsilon (n+2\mu) \\
  \end{array}
\right).$$
On the other hand, it follows from the definition  of {\it semisimple} that $A$ can be diagonalized. Thus  it is necessary that  $b_0$ and $c_0(\la n-2\mu)$ are equal to $0$.

By Lemma \ref{lemma5}, one can write
\begin{eqnarray}\label{jj}
\left\{\begin{array}{llll}
\theta^{-1}\sigma(L_n)=\epsilon \bar{\al}^n L_{\epsilon n}+\sum\limits_{i'\in \bar{S}}\bar{b}_{n_{i'}}Y_{i'}+\sum\limits_{j'\in \bar{S'}}\bar{c}_{n_{j'}}M_{j'},\vs{6pt}\\
\theta^{-1}\sigma(Y_n)=\bar{\al}^n\bar{\beta}Y_{\epsilon n}+\sum\limits_{k\in \bar{S}''}\bar{e}_{n_k}M_k,\vs{6pt}\\
\theta^{-1}\sigma(M_n)=\epsilon \bar{\al}^n\bar{\beta}^2M_{\epsilon n},
\end{array}\right.
\end{eqnarray}
for some $\bar{\al}$, $\bar{\beta}\in\C^*$, $\bar{b}_{n_{i'}}$, $\bar{c}_{n_{j'}}$ and $\bar{e}_{n_k}\in\C$, $\epsilon \in\{\pm1\}$ with $\epsilon=1$ for $\mu\notin\frac{1}{2}\Z$, $\bar{S}$, $\bar{S'}$, $\bar{S''}\subseteq\Z$.

Applying $\theta^{-1}\sigma$ to $(n+\mu)Y_n=[L_0,Y_n]$, we obtain
\begin{eqnarray*}
&&(n+\mu)\Big(\bar{\al}^n\bar{\beta}Y_{\epsilon n}+\mbox{$\sum\limits_{k\in \bar{S}''}$}\bar{e}_{n_k}M_k\Big)\\
\!\!\!\!&\!\!\!\!\!\!\!\!\!=\!\!\!\!\!\!\!\!\!&[\epsilon L_0+c_0M_0,\bar{\al}^n\bar{\beta}Y_{\epsilon n)}+\mbox{$\sum\limits_{k\in \bar{S}''}$}\bar{e}_{n_k}M_k]=(n+\mu)\bar{\al}^n\bar{\beta}Y_{\epsilon n}+\mbox{$\sum\limits_{k\in \bar{S}''}$}\epsilon\bar{e}_{n_k}(k+2\mu)M_k,
\end{eqnarray*}
 which gives
\begin{eqnarray*}
(n+\mu)\mbox{$\sum\limits_{k\in \bar{S}''}$}\bar{e}_{n_k}M_k=\mbox{$\sum\limits_{k\in \bar{S}''}$}\epsilon\bar{e}_{n_k}(k+2\mu)M_k.
\end{eqnarray*}
Comparing the coefficients of $M_k$, one has
\begin{equation}\label{0052}
\bar{e}_{n_k}\big(k+2\mu-\epsilon(n+\mu)\big)=0,\ \ \forall\ k\in \bar{S''}.
\end{equation}

Applying $\theta^{-1}\sigma$ to $nL_n=[L_0,L_n]$, we obtain
\begin{eqnarray}\label{0058}
\nonumber&&n\Big(\epsilon \bar{\al}^n L_{\epsilon n}+\mbox{$\sum\limits_{i'\in \bar{S}}$}\bar{b}_{n_{i'}}Y_{i'}+\mbox{$\sum\limits_{j'\in \bar{S'}}$}\bar{c}_{n_{j'}}M_{j'}\Big)\\
&=&[\epsilon L_0+c_0M_0,\epsilon \bar{\al}^n L_{\epsilon n}+\mbox{$\sum\limits_{i'\in \bar{S}}$}\bar{b}_{n_{i'}}Y_{i'}+\mbox{$\sum\limits_{j'\in \bar{S'}}$}\bar{c}_{n_{j'}}M_{j'}]\\
\nonumber&=&\epsilon n\bar{\al}^nL_{\epsilon n}+\mbox{$\sum\limits_{i'\in \bar{S}}$}\epsilon \bar{b}_{n_{i'}}(i'+\mu)Y_{i'}+\mbox{$\sum\limits_{j'\in \bar{S'}}$}\epsilon \bar{c}_{n_{j'}}(j'+2\mu)M_{j'}+\epsilon\bar{\al}^nc_0(\la\epsilon n-2\mu)M_{\epsilon n}.
\end{eqnarray}
Using $c_0(\la n-2\mu)=0$ and $\epsilon=1$ for $\mu\notin\frac{1}{2}\Z$, one can deduce $\epsilon\bar{\al}^nc_0(\la\epsilon n-2\mu)=0$. Using this in (\ref{0058}), one has
\begin{eqnarray*}
\nonumber&&n\Big(\mbox{$\sum\limits_{i'\in \bar{S}}$}\bar{b}_{n_{i'}}Y_{i'}+\mbox{$\sum\limits_{j'\in \bar{S'}}$}\bar{c}_{n_{j'}}M_{j'}\Big)=\mbox{$\sum\limits_{i'\in \bar{S}}$}\epsilon \bar{b}_{n_{i'}}(i'+\mu)Y_{i'}+\mbox{$\sum\limits_{j'\in \bar{S'}}$}\epsilon \bar{c}_{n_{j'}}(j'+2\mu)M_{j'}.
\end{eqnarray*}
Comparing the coefficients of $Y_{i'}$ and $M_{j'}$ respectively, we have
\begin{eqnarray}\label{0053}
\bar{b}_{n_{i'}}(i'-\epsilon n+\mu)=0,\ \forall\,i'\in \bar{S},\ \ \bar{c}_{n_{j'}}(j'-\epsilon n+2\mu)=0,\ \ \forall\,j'\in \bar{S'}.
\end{eqnarray}
\setcounter{case}{0}
\begin{case}
$\mu\notin\frac{1}{2}\Z$.
\end{case}
Since $\mu\notin\frac{1}{2}\Z$, we have $\epsilon=1$. Thus, by  (\ref{0052}) and (\ref{0053}), we obtain  $\bar{b}_{n_{i'}}$, $\bar{c}_{n_{j'}}$ and $\bar{e}_{n_k}$, are all equal to zero for all $i'\in \bar{S}$, $j'\in \bar{S'}$ and $k\in \bar{S''}$. Using this in (\ref{jj}), we obtain
\begin{eqnarray*}
\theta^{-1}\sigma(L_n)=\bar{\al}^n L_{n},\ \ \theta^{-1}\sigma(Y_n)=\bar{\al}^n\bar{\beta}Y_{n},\ \ \theta^{-1}\sigma(M_n)=\bar{\al}^n\bar{\beta}^2M_{n}.
\end{eqnarray*}
Let ${\varphi}_{\bar{\al},\bar{\beta}}$ be the automorphism of $\mathscr{L}$ as that given in Lemma \ref{auto} (ii). Then $\sigma=\theta\varphi_{\bar{\al},\bar{\beta}}$.
\begin{case}
$\mu=0$.
\end{case}
Since $\mu=0$, we know that ${i'}$, ${j'}$ and $k$ are all equal to $\epsilon n$ for all  $i'\in \bar{S}$, $j'\in \bar{S'}$ and $k\in \bar{S''}$ in (\ref{0052}) and (\ref{0053}). Using this in (\ref{jj}), we obtain
\begin{eqnarray}\label{0070}
\left\{\begin{array}{lll}
\theta^{-1}\sigma(L_n)=\epsilon \bar{\al}^nL_{\epsilon n}+\bar{b}_nY_{\epsilon n}+\bar{c}_nM_{\epsilon n},\vs{6pt}\\
\theta^{-1}\sigma(Y_n)=\bar{\al}^n\bar{\beta}Y_{\epsilon n}+\bar{e}_nM_{\epsilon n},\vs{6pt}\\
\theta^{-1}\sigma(M_n)=\epsilon \bar{\al}^n\bar{\beta}^2M_{\epsilon n},
\end{array}\right.
\end{eqnarray}
for some $\bar{b}_n$, $\bar{c}_n$ and $\bar{e}_n\in\C$. Let $\bar{\theta}$ be the automorphism of $\mathscr{L}$ defined by
\begin{equation*}
\bar{\theta}(X_n)=\bar{\al}^nX_n\ \ \ {\rm for}\ X\in\{L,Y,W\}.
\end{equation*}
Then $\bar{\theta}$ is an inner one. Set $$\bar{f}(n)=\bar{b}_n(\bar{\al}^n\bar{\beta})^{-1},\ \ f(n)=\bar{c}_n(\epsilon\bar{\al}^n\bar{\beta}^2)^{-1},\ \ g(n)=\bar{e}_n(\epsilon\bar{\al}^n\bar{\beta}^2)^{-1}.$$
Define $\varphi_\epsilon$ and
${\varphi}_{1,\bar{\beta}}$ as those given in Lemma \ref{auto} (i) and (ii) respectively. Then we can write  (\ref{0070}) as follows
\begin{eqnarray*}
\phi(L_n)=L_{n}+\bar{f}(n)Y_n+f(n)M_n,\ \ \phi(Y_n)=Y_{n}+g(n)M_n,\ \
\phi(M_n)=M_{n},
\end{eqnarray*}
where $\phi=(\varphi_\epsilon)^{-1}({\varphi}_{1,\bar{\beta}})^{-1}(\bar{\theta})^{-1}\theta^{-1}\sigma$.
\setcounter{subcase}{0}
\begin{subcase}
$\la=-2$.
\end{subcase}
By Lemma \ref{lemma9}, there exist some $c_1$, $\bar{c}_1$ and $b_1\in\C$ such that
\begin{eqnarray*}
&&\phi(L_n)=L_{n}+b_1nY_n+\big(\frac{c_1}{6}(n^3-n)-\frac{\bar{c}_1}{3} (n^3-4n)+\frac{b_1^2}{3}n(n-2)(n-1)\big)M_n,\\
&&\phi(Y_n)=Y_{n}-2b_1nM_n,\ \ \phi(M_n)=M_{n}.
\end{eqnarray*}
Set $\theta_1={\rm exp}\big(-2b_1{\rm ad}Y_0+ \frac{1}{12}(c_1-8\bar{c}_1-4b_1^2){\rm ad}M_0\big)$. Then
\begin{eqnarray*}
(\theta_1)^{-1}\phi(L_n)=L_{n}+\frac{c_1-2\bar{c}_1+2b_1^2}{6}n^3M_n,\ \ (\theta_1)^{-1}\phi(Y_n)=Y_{n},\ \ (\theta_1)^{-1}\phi(M_n)=M_{n}.
\end{eqnarray*}
Set $\al_1=\frac{1}{6}(c_1-2\bar{c}_1+2b_1^2)$ and define $\phi_{\al_1}$ as in Lemma \ref{auto} (iii). Then $(\phi_{\al_1})^{-1}(\theta_1)^{-1}\phi=Id$.
\begin{subcase}
$\la=-1$.
\end{subcase}
By Lemma \ref{lemma9}, there exist some $c_2$, $\bar{c}_2$ and $b_2\in\C$ such that
\begin{eqnarray*}
\phi(L_n)=L_{n}+\big(\frac{c_2}{2}(n^2-n)-\bar{c}_2(n^2-2n)\big)M_n,\ \ \phi(Y_n)=Y_{n}+b_2nM_n,\ \ \phi(M_n)=M_{n}.
\end{eqnarray*}
Define $\psi_{b_2}$ that given as in Lemma \ref{auto} (iv). Then
\begin{eqnarray*}
&&(\psi_{b_2})^{-1}\phi(L_n)=L_{n}+\big(\frac{c_2}{2}(n^2-n)-\bar{c}_2(n^2-2n)\big)M_n,\\ &&(\psi_{b_2})^{-1}\phi(Y_n)=Y_{n},\ \ (\psi_{b_2})^{-1}\phi(M_n)=M_{n}.
\end{eqnarray*}
Set $\theta_2={\rm exp}\ \big(\frac{1}{2}(c_2-4\bar{c}_2){\,\rm ad}M_0\big)$. Then
\begin{eqnarray*}
&&(\theta_2)^{-1}(\psi_{b_2})^{-1}\phi(L_n)=L_{n}+\frac{c_2-2\bar{c}_2}{2}n^2M_n,\\ &&(\theta_2)^{-1}(\psi_{b_2})^{-1}\phi(Y_n)=Y_{n},\ \ (\theta_2)^{-1}(\psi_{b_2})^{-1}\phi(M_n)=M_{n}.
\end{eqnarray*}
Set $\al_2=\frac{1}{2}(c_2-2\bar{c}_2)$, let $\phi_{\al_2}$ be as that given in Lemma \ref{auto}(iii). Then $(\phi_{\al_2})^{-1}(\theta_2)^{-1}(\psi_{b_2})^{-1}\phi=Id$.
\begin{subcase}
$\la=1$.
\end{subcase}
By Lemma \ref{lemma9}, there exist some $e$, $c_4$, $\bar{c}_4$ and $b_4\in\C$ such that
\begin{eqnarray*}
&&\phi(L_n)=L_{n}+b_4nY_n+\big(c_4n+\frac{b_4^2}{2}n(n-1)\big)M_n,\\ &&\phi(Y_n)=Y_{n}+(b_4n+e)M_n,\ \ \phi(M_n)=M_{n}.
\end{eqnarray*}
Set $\theta_4={\rm exp}\,\big(b_4{\rm ad}Y_0- \frac{1}{2}(b_4^2-2c_4){\,\rm ad}M_0\big)$. Then
\begin{eqnarray*}
&&(\theta_4)^{-1}\phi(L_n)=L_{n},\ \ (\theta_4)^{-1}\phi(Y_n)=Y_{n}+eM_n,\ \ (\theta_4)^{-1}\phi(M_n)=M_{n}.
\end{eqnarray*}
Let $\psi_e$ be as that given in Lemma {\ref{auto}} (iv). Then $(\psi_e)^{-1}(\theta_4)^{-1}\phi=Id$.
\begin{subcase}
$\la\notin\{-2,0,\pm1\}$.
\end{subcase}
By Lemma \ref{lemma9}, there exist some $c_5$, $\bar{c}_5$ and $b_5\in\C$ such that
\begin{eqnarray*}
&&\phi(L_n)=L_{n}+b_5nY_n+\big(\bar{c}_5n+\frac{b_5^2}{\la+1}n(n-1)\big)M_n,\\ &&\phi(Y_n)=Y_{n}+\frac{2b_5}{\la+1}nM_n,\ \ \phi(M_n)=M_{n}.
\end{eqnarray*}
Set $\theta_5={\rm exp}\,\big(2b_5(\la+1)^{-1}{\rm ad}Y_0+ (\bar{c}_5(\la+1)-b_5^2)(\la^2+\la)^{-1}{\rm ad}M_0\big)$. Then
$(\theta_5)^{-1}\phi=Id$.

Recall that $\phi=(\varphi_\epsilon)^{-1}({\varphi}_{1,\bar{\beta}})^{-1}(\bar{\theta})^{-1}\theta^{-1}\sigma$ . Since $Inn\mathscr{L}$ is the normal subgroup of $Aut\mathscr{L}$, then there exist $\bar{\theta}_i\in Inn\mathscr{L}(i=1,2,\cdots,5)$ such that
\begin{eqnarray*}
\sigma=\left\{\begin{array}{llllllll}%
\bar{\theta}_1{\varphi}_{1,\bar{\beta}}\varphi_\epsilon\phi_{\al_1}&{\rm if}\ \ \la=-2,\vs{6pt}\\
\bar{\theta}_2{\varphi}_{1,\bar{\beta}}\varphi_\epsilon\phi_{\al_2}\psi_{b_2}&{\rm if}\ \ \la=-1,\vs{6pt}\\
\bar{\theta}_4{\varphi}_{1,\bar{\beta}}\varphi_\epsilon\psi_{e}&{\rm if}\ \ \la=1,\vs{6pt}\\
\bar{\theta}_5{\varphi}_{1,\bar{\beta}}\varphi_\epsilon&{\rm if}\ \ \la\notin\{-2,0,\pm1\}.
\end{array}\right.
\end{eqnarray*}

Obviously, $Inn\mathscr{L}$ satisfies (\ref{inn}). Hence the theorem follows from Lemma {\ref{lemma9}}.\QED

\end{document}